\documentclass[11pt,twoside]{article}
\usepackage{amsfonts}
\usepackage{amsmath,amssymb}
\usepackage{amsthm}
\usepackage[mathscr]{euscript}
\usepackage{mathrsfs}
\usepackage{indentfirst}
\usepackage{booktabs,multirow,bigstrut,makecell}
\usepackage{color}\usepackage{enumitem}

\numberwithin{equation}{section}

\newcommand{\lbl}[1]{\label{#1}}

\newtheorem{theo}{Theorem}[section]

\newtheorem{rem}{Remark}[section]

\allowdisplaybreaks

\newcommand{\be}{\begin{equation}}
\newcommand{\ee}{\end{equation}}
\newcommand\bes{\begin{eqnarray}} \newcommand\ees{\end{eqnarray}}
\newcommand\beslf{\begin{eqnarray}\left\{\begin{array}{lll}}
\newcommand\eesrr{\end{array}\right.\end{eqnarray}}
\newcommand{\bess}{\begin{eqnarray*}}
\newcommand{\eess}{\end{eqnarray*}}

\newcommand\dd{\displaystyle}
\newcommand\df{\dd\frac}

\newcommand\yy{\infty}

\newcommand\R{\mathbb{R}}
\newcommand\up{\overline}

\newcommand\n{\nonumber}
\newcommand\kk{\left}
\newcommand\ol{\overline}
\newcommand\rr{\right}

\newcommand\sk{\smallskip}
\newcommand\mk{\medskip}

\newcommand{\bbbb}{\left\{\begin{array}{lllll}}
\newcommand{\nnnn}{\end{array}\right.}

\begin{document}\thispagestyle{empty}
\setlength{\baselineskip}{16pt}

\begin{center}{\Large\bf Existence and uniqueness of solution of free}\\[2mm]
{\Large\bf boundary problems with partially degenerate diffusion}\footnote{This work was
supported by NSFC Grants 11771110, 11971128}\\[4mm]
 {\Large  Siyu Liu, \ Mingxin Wang\footnote{Corresponding author. {\sl E-mail}: mxwang@hit.edu.cn}}\\[1.5mm]
School of Mathematics, Harbin Institute of Technology, Harbin 150001, PR China
\end{center}

\begin{quote}
\noindent{\bf Abstract.}
In this paper, we mainly introduce a general method to study the existence and uniqueness of solution of free boundary problems with partially degenerate diffusion.


\end{quote}

 \section{Introduction and main result}
  \setcounter{equation}{0}{\setlength\arraycolsep{2pt}

In recent years, the following free boundary problem with a partially degenerate diffusion
  \be\label{1.1}
\begin{cases}
u_t=f_1(t,x,u,v),&t>0,\ \ g(t)<x<h(t),\\[1mm]
v_t=dv_{xx}+f_2(t,x,u,v),&t>0,\ \ g(t)<x<h(t),\\[1mm]
u(t,x)=v(t,x)=0,&t\ge0, \ \ x=g(t),\ h(t),
\\[1mm]
g'(t)=-\mu v_x(t,g(t)),\ \ h'(t)=-\beta v_x(t,h(t)),&t\ge0,\\[1mm]
u(0,x)=u_0(x),\ \ v(0,x)=v_0(x),&-h_0\le x\le h_0,\\[1mm]
h(0)=-g(0)=h_0,
 \end{cases}\ee
has been studied by some authors to describe the nature of spreading and vanishing of multiple species, where $d, \mu,\beta$ and $h_0$ are positive constants. In the problem \eqref{1.1}, the diffusion of species $u$ is relatively faster than that of species $v$, or species $v$ has no diffusion, so the diffusion of $v$ is omitted. Wang and Cao (\cite{wc15}) studied the case $f_i(t,x,u,v)=f_i(u,v)$ and $\beta=\mu$, where $(f_1(u,v),f_2(u,v))$ has a cooperative structure and is controlled by a linear system. Ahn et al. (\cite{Ahn}) investigated a man-environment-man epidemic model:
$f_1(t,x,u,v)=G(v)-au$, $f_2(t,x,u,v)=bu-cv$ and $\beta=\mu$. Tarboush et al. (\cite{tl17}) discussed a West Nile virus model: $f_1(t,x,u,v)=r_1(a-u)v-bu$, $f_2(t,x,u,v)=r_2(b-v)u-cv$  and $\beta=\mu$. In the study of the local existence and uniqueness of solution, they used the different methods. In \cite{Ahn}, the function $G$ satisfies

\vskip 4pt $\bullet$ \ $G\in C^1([0,\yy))$, $G(0)=0$, $G'(v)>0$, $\frac{G(v)}v$ is decreasing and $\lim\limits_{v\to\yy}\frac{G(v)}v<{ac}/b$.

\vskip 4pt In \eqref{1.1}, the curves $x=g(t)$ and $x=h(t)$ are the free boundaries to be determined together with $u(t,x)$ and $v(t,x)$.

The main aim of this paper is to give another rigorous proof of existence and uniqueness of solution. Denote $C^{1-}(\Omega)$ be the Lipschitz continuous functions space. We assume that the initial functions $u_0$ and $v_0$ satisfy
\begin{itemize}
  \item $(u_0,v_0)\in C^{1-}([-h_0,h_0])\times W_p^2((-h_0,h_0))$ with $p>3$, $u_0(\pm h_0)=v_0(\pm h_0)=0$, $u_0,v_0>0$ in $(-h_0,h_0)$, and $v'_0(h_0)<0, v'_0(-h_0)>0$,
\end{itemize}
and denote by $L_0$ the Lipschitz constant of $u_0$ in $x$.

It is assumed that $(f_1,f_2)$ satisfies
\begin{description}
  \item[{\bf (I)}]  $f_1(t,x,0,v)\ge 0$ for all $v\ge 0$ and $f_2(t,x,u,0)\ge 0$ for all
  $u\ge 0$. For any given $\tau$, $l$, $k_1$, $k_2>0$, $f_i(\cdot,0,0)\in L^\infty((0,\tau)\times(-l,l))$ and there exists a constant $L(\tau,l, k_1, k_2)>0$ such that
 \bess
|f_i(t,x,u_1,v_1)-f_i(t,x,u_2,v_2)|\le L(\tau,l, k_1, k_2)(|u_1-u_2|+|v_1-v_2|),\ i=1,2
\eess
for all $t\in[0,\tau]$, $x\in[-l,l]$, $u_1,u_2\in[0,k_1]$, $v_1,v_2\in[0,k_2]$;
  \item[{\bf (II)}] $f_{i}$ is locally Lipschitz continuous in $x\in\R$, i.e., for any given any
given $\tau$, $l$, $k_1$, $k_2>0$, there exists a constant $L^*(\tau,l, k_1, k_2)>0$ such that
\bess
|f_i(t,x,u,v)-f_i(t,y,u,v)|\le L^*(\tau,l, k_1, k_2)|x-y|,\ \ i=1,2
\eess
for all $t\in[0,\tau]$, $x,\,y\in[-l,l]$, $u\in[0,k_1]$, $v\in[0,k_2]$.
\end{description}

It is easy to notice that the condition {\bf (I)} implies $f_i\in L^\infty((0,\tau)\times(-l,l)\times(0,k_1)\times(0,k_2))$ for any given $\tau$, $l$, $k_1$, $k_2>0$.

The result concerns with the local existence and uniqueness.

\begin{theo}\label{th1.1}\, Under the above assumptions, there exists a $T>0$ such that the problem \eqref{1.1} has a unique solution $(u,v,g,h)$ which is defined on $[0,T]$. Moreover,
 \bess
 &g,h\in C^{1+\frac \alpha 2}([0,T]),\ \ g'(t)<0,\ h'(t)>0 \ \ \mbox{in} \ \ [0,T],&\\[1mm]
 &u\in C^{1,1-}(\overline D^T_{g,h}), \ v\in W_p^{1,2}(D^T_{g,h})\cap C^{1+\alpha,\frac{1+\alpha}{2}}(\overline D^T_{g,h}),\ \ u,\ v>0 \ \ \mbox{in} \ \
 D^T_{g,h},&\eess
where
 \[D^T_{g,h}=\{0<t\le T,\ g(t)<x<h(t)\},\]
$u\in C^{1,1-}(\overline D^T_{g,h})$ means that $u$ is differentiable continuously in $t\in[0,T]$ and is Lipschitz continuous in $x\in[g(t),h(t)]$.
 \end{theo}

When $f_1, f_2$ do not depend on $(t,x)$, i.e., $f_1(t,x,u,v)=f_1(u,v)$, $f_2(t,x,u,v)=f_2(u,v)$, we have the following global existence results.

\begin{theo}\label{th2.1}\, Let $(f_1,f_2)$ be qusi-monotone increasing for $u,v\ge 0$. If the initial value problem
 \bess\left\{\begin{array}{ll}
 \phi'(t)=f_1(\phi,\psi), \ \ \psi'(t)=f_2(\phi,\psi), \ \ t>0,\\[2mm]
 \phi(0)=\dd\max_{[-h_0,h_0]} u_0>0, \ \ \psi(0)=\max_{[-h_0,h_0]} v_0>0
 \end{array}\right.
 \eess
has a global solution $(\phi,\psi)$, then the unique solution $(u,v,g,h)$ of \eqref{1.1} also exists globally.
\end{theo}

\begin{theo}\label{th2.2}\, Assume that there exists $k_0>0$ such that $f_1(u,v)<0$ for all $u>k_0$, $v\ge 0$, and for the given $\eta>0$, there exists $\Theta(\eta)>0$ such that $f_2(u,v)<0$ for $0\le u\le \eta$, $v\ge\Theta(\eta)$, then the unique solution $(u,v,g,h)$ of \eqref{1.1} exists globally.
\end{theo}

\begin{rem} {\rm (i)} Conditions $u(t,g(t))=u(t,h(t))=0$ in \eqref{1.1} look like
boundary conditions of $u$, but they do actually play the roles of initial conditions of $u$ at points $x=g(t)$ and $x=h(t)$, respectively.

{\rm(ii)} Our conclusions are applicable to the models investigated in {\rm\cite{Ahn, tl17, wc15}},  and assert that the solution exists globally {\rm(}using Theorems {\rm\ref{th1.1}} and {\rm\ref{th2.1}} for the models in {\rm\cite{Ahn,wc15}}, and Theorems {\rm\ref{th1.1}} and {\rm\ref{th2.2}} for the models in {\rm\cite{tl17})}.
\end{rem}

\section{Proofs of Theorems \ref{th1.1}-\ref{th2.2}}

{\bf Proof of Theorem \ref{th1.1}}. The proof is divided into several steps.

{\it Step 1:} For $0<T<\yy$, set
 \[A=\max_{[-h_0,h_0]} u_0+1, \ \ B=\max_{[-h_0,h_0]} v_0+1, \ \ \Pi_T=[0,T]\times[-2h_0, 2h_0],\ \ \ \Delta_T=[0,T]\times[-1,1],\]
and denote $L_1=L(1,2h_0, A, B)$,
 \[{\cal A}=\{d, h_0, \mu, \beta, A, B, \|v_0\|_{W_p^2((-h_0,h_0))}, v_0'(\pm h_0), \|f_2\|_{L^\yy(\Pi_1\times(0,A)\times(0,B))}, L_1\}.\]
We say $u\in C^{1-}_{x}(\Pi_T)$ if there is a constant $L(u,T)$ such that
 \[|u(t,x_1)-u(t,x_2)|\le L(u,T)|x_1-x_2|, \ \ \forall \ x_1,x_2\in[-2h_0, 2h_0], \ t\in[0,T].\]
Define
  \bess
  \mathbb X_{u_0}^T&=&\left\{\phi\in C(\Pi_T): \phi(0,x)=u_0(x),\ 0\le\phi\le
  A\right\}.
  \eess
Chosen $u\in\mathbb X_{u_0}^1\cap C^{1-}_{x}(\Pi_1)$ and consider the following problem
 \bes
\left\{\begin{aligned}
&v_t=d v_{xx}+f_2(t,x,u(t,x),v), &&0<t\le 1,~g(t)<x<h(t),\\
&v(t,g(t))=v(t,h(t))=0, &&0\le t\le 1,\\
&g'(t)=-\mu v_x(t,g(t)), \ h'(t)=-\beta v_x(t,h(t)),&&0\le t\le 1,\\
&v(0,x)=v_0(x),\ h(0)=-g(0)=h_{0}>0, &&|x|\le h_{0}.
\end{aligned}\right.
 \label{2.1}
 \ees
Due to the properties of $f_2$ and $u$, using the similar arguments in the proof of \cite[Theorem 1.1]{Wdcds19} we can show that there exists $0<T_0\ll 1$ such that \eqref{2.1} has a unique solution $(v,g,h)$ and satisfies
 \[g, h\in C^{1+\frac\alpha 2}([0,T_0]), \ \ \ v\in W^{1,2}_p(D^{T_0}_{g,h})\cap  C^{\frac{1+\alpha}{2},\,1+\alpha}(\ol D^{T_0}_{g,h})\]
with $0<\alpha<1-3/p$, and
  \bes\label{2.2}
  \begin{cases}
   &\|g,h\|_{C^{1+\frac \alpha 2}([0,T_0])},\ \|v\|_{W^{1,2}_p(D^{T_0}_{g,h})}\leq K, \ \
    0<v\le B \ \ {\rm in}\ \ D^{T_0}_{g,h},\ \ \ \\[2mm]
  &0<-g'(t), \ h'(t)\le K, \ \ \ |g(t)|, \ h(t)\le 2h_0\ \  \mbox{on}\ \ [0,T_0],
\end{cases}
\ees
where $T_0$ and $K$ depend only on ${\cal A}$, $\alpha$ and the Lipschitz
constant $L(u,1)$ of $u$.

Define $\tilde u_0(x)=u_0(x)$ when $|x|\le h_0$, and $\tilde u_0(x)=0$ when $|x|>h_0$. Then
$\tilde u_0\in C^{1-}([g(T_0), h(T_0)])$ since $u_0\in C^{1-}([-h_0,h_0])$.
For the functions $g(t)$, $h(t)$ obtained above, it is easy to see that the inverse functions
$g^{-1}(x)$ and $h^{-1}(x)$ exist for $x\in[g(T_0), h(T_0)]$. We set
  \bes
t_x=\left\{\begin{array}{ll}
g^{-1}(x) \ \ &\mbox{if }\ x\in[g(T_0),-h_0),\\
0 &\mbox{if } \ |x|\le h_0,\\
h^{-1}(x)\ \ &\mbox{if } \ x\in(h_0,h(T_0)],
 \end{array}\right.
 \lbl{2.3}\ees
 which is Lipschitz continuous in $x$. For the function $v(t,x)$ obtained above, and every $g(T_0)<x<h(T_0)$, we consider the following problem
 \be\label{2.4}
 \begin{cases}
\tilde u_t=f_1(t,x,\tilde u,v(t,x)),&t_x<t\le T_0,\\[1mm]
\tilde u(t_x,x)=\tilde u_0(x).
 \end{cases}\ee
By the standard theory for ODE we can see that there exists $0<T<T_0$, which depends on $L_1$, $A$, $B$ and $K$, such that for all $g(T)\le x\le h(T)$, $\tilde u(t,x)$ is defined on $[t_x,T]$ and so $\tilde{u}$ is defined on $\ol D^{T}_{g,h}$. Moreover, as the function of $(t,x)$, we assert that $\tilde u\in C^{1,1-}(\ol D^{T}_{g,h})$, of which the Lipschitz constant will be calculated in the next step, and $\tilde u(t,g(t))=\tilde u(t,h(t))=0$ for $0\le t\le T$. Make the zero extension of $\tilde u$ to $[0,t_x]$ for every $g(T)\le x\le h(T)$. Then $\tilde u\in C^{1,1-}([0,T]\times[g(T), h(T)])$.

\mk{\it Step 2:} The estimate of the Lipschitz constant of $\tilde u$ in $x$. As $g'(t)<0, h'(t)>0$  on $[0,T_0]$, there is a $\sigma>0$ such that
 \bes
 |g'(t)|\ge\sigma,  \ \ |h'(t)|\ge \sigma,  \ \ \forall \ t\in[0,T_0].
 \lbl{2.5}\ees
Set $F_1(s,x)=f_1(s,x,\tilde{u}(s,x),v(s,x))$. It follows from the first equation of \eqref{2.4} that,  for $t_x<\tau\le T$,
 $$ \tilde{u}(\tau,x)=\tilde{u}(t_x,x)+\int_{t_x}^\tau F_1(s,x){\rm d}s.$$
For the given $(t,x_1), (t,x_2)\in\Pi_T$. We divide the arguments into several cases

\mk{\bf Case 1}: $(t,x_1), (t,x_2)\in\ol D^{T}_{g,h}$ with $-h_0\le x_2<x_1$. Then $t\ge t_{x_1}\ge t_{x_2}\ge 0$. Thus we have, for any $t_{x_1}\le\tau\le t$,
  \begin{align*}
  |\tilde{u}(\tau,x_1)-\tilde{u}(\tau,x_2)|\le |\tilde{u}(t_{x_1},x_1)-\tilde{u}(t_{x_2},x_2)|+\int_{t_{x_2}}^{t_{x_1}}|F_1(s,x_2)|{\rm d}s+\int_{t_{x_1}}^\tau|F_1(s,x_1)-F_1(s,x_2)|{\rm d}s.
  \end{align*}
By use of the conditions {\bf (I)} and {\bf (II)}, it is easy to derive that
\begin{align*}
|F_1(s,x_1)-F_1(s,x_2)|\le L_1\big(|\tilde{u}(s,x_1)-\tilde{u}(s,x_2)|+|v(s,x_1)-v(s,x_2)|\big)+L^*_1|x_1-x_2|,
\end{align*}
where $L_1=L(1,2h_0, A, B)$, $L_1^*=L^*(1,2h_0, A, B)$. It yields,
\begin{align}
|\tilde{u}(\tau,x_1)-\tilde{u}(\tau,x_2)|\le& T L_1 \|\tilde u(\cdot,x_1)-\tilde u(\cdot,x_2)\|_{C([t_{x_1},t])}+|\tilde{u}(t_{x_1},x_1)-\tilde{u}(t_{x_2},x_2)|\n\\[1mm]
&+L^*_1|x_1-x_2|+C_1|t_{x_1}-t_{x_2}|+L_1\int_{t_{x_1}}^\tau
|v(s,x_1)-v(s,x_2)|{\rm d}s\label{2.6}
\end{align}
as $\tau\le T\le 1$, where $C_1=\|f_1\|_{L^\yy(\Pi_1\times(0,A)\times(0,B))}$.

Noticing $\|v\|_{W^{1,2}_p(D^{T_0}_{g,h})}\le K$, we have $\|v_x\|_{L^{\yy}(\up{D}^{T}_{g,h})}\le C_2$ by the embedding theorem. Thus
 \bess
\int_{t_{x_1}}^\tau|v(s,x_1)-v(s,x_2)|{\rm d}s\le T\|v_x\|_{L^{\yy}(\up{D}^{T}_{g,h})}|x_1-x_2|
\le TC_2|x_1-x_2|\le C_2|x_1-x_2|
 \eess
as $\tau\le T\le 1$.

If $t_{x_1}>0$, $t_{x_2}>0$, then $\tilde{u}(t_{x_1},x_1)=\tilde{u}(t_{x_2},x_2)=0$ and
 \bess
|t_{x_1}-t_{x_2}|=|h^{-1}(x_1)-h^{-1}(x_2)|\le\|(h^{-1})'\|_{L^{\yy}([0,T])}|x_1-x_2|
\le \sigma^{-1}|x_1-x_2|.
 \eess

If $t_{x_1}>0$, $t_{x_2}=0$, then $x_2\in[-h_0,h_0]$, $x_1>h_0$, $\tilde{u}(t_{x_1},x_1)=0$. Let $L_0$ be the Lipschitz constant of $u_0$ in $x$. It then follows that
 \bess
&|t_{x_1}-t_{x_2}|=|h^{-1}(x_1)-0|=|h^{-1}(x_1)-h^{-1}(h_0)|\le \sigma^{-1}|x_1-h_0|\le \sigma^{-1}|x_1-x_2|,&\\[1mm]
&|\tilde{u}(t_{x_1},x_1)-\tilde{u}(t_{x_2},x_2)|=|0-u_0(x_2)|=|u_0(h_0)-u_0(x_2)|\le L_0|h_0-x_2|\le L_0|x_1-x_2|.&
 \eess

If $t_{x_1}=t_{x_2}=0$, i.e., $x_1,x_2\in[-h_0,h_0]$, then
\bess
|\tilde{u}(t_{x_1},x_1)-\tilde{u}(t_{x_2},x_2)|=|u_0(x_1)-u_0(x_2)|\le L_0|x_2-x_1|.
\eess

Substituting these estimates into \eqref{2.6}, we have
\begin{align*}
|\tilde{u}(\tau,x_1)-\tilde{u}(\tau,x_2)|\le TL_1\|\tilde{u}(\cdot,x_1)-\tilde{u}(\cdot,x_2)\|_{C([t_{x_1},t])}+(L_0+L^*_1
+C_1\sigma^{-1}+C_2L_1)|x_1-x_2|.
\end{align*}
Take the maximum of $|\tilde{u}(\tau,x_1)-\tilde{u}(\tau,x_2)|$ in $[t_{x_1},t]$ it yields
\bess
 &&\|\tilde{u}(\cdot,x_1)-\tilde{u}(\cdot,x_2)\|_{C([t_{x_1},t])}\\[1mm]
 &\le& TL_1\|\tilde{u}(\cdot,x_1)-\tilde{u}(\cdot,x_2)\|_{C([t_{x_1},t])}+(L_0+L^*_1
+C_1\sigma^{-1}+C_2L_1)|x_1-x_2|.
\eess
Set $M=2(L_0+L_1^*+C_1\sigma^{-1}+C_2L_1)$. Then
 \bes
|\tilde{u}(t,x_1)-\tilde{u}(t,x_2)|\le\|\tilde{u}(\cdot,x_1)-\tilde{u}(\cdot,x_2)\|_{C([t_{x_1},t])}\le M|x_1-x_2|\lbl{2.7}
 \ees
provided that $0<T\le \min\{1,\frac{1}{2L_1}\}$.\vskip 4pt

\mk{\bf Case 2}: $(t,x_1), (t,x_2)\in\ol D^{T}_{g,h}$ with $x_2<x_1\le h_0$. Similar to the above, \eqref{2.7} holds.

\mk{\bf Case 3}: $(t,x_1), (t,x_2)\in\ol D^{T}_{g,h}$ with $x_2<-h_0<h_0<x_1$. Then
 \bes
 |\tilde{u}(t,x_1)-\tilde{u}(t,x_2)|\le|\tilde{u}(t,x_1)-\tilde{u}(t,h_0)|
 +|\tilde{u}(t,h_0)-\tilde{u}(t,x_2)|\le 2M|x_1-x_2|.
 \label{2.8}\ees

{\bf Case 4}: $(t,x_1), (t,x_2)\not\in\ol D^{T}_{g,h}$. Then $\tilde{u}(t,x_1)=\tilde{u}(t,x_2)=0$.

\mk{\bf Case 5}: $(t,x_1)\not\in\ol D^{T}_{g,h}$, $(t,x_2)\in\ol D^{T}_{g,h}$. We may assume that $x_1> h(t)$. Thus $x_2\le h(t)$, $\tilde{u}(t,x_1)=\tilde{u}(t,h(t))=0$, and
 \[|\tilde{u}(t,x_1)-\tilde{u}(t,x_2)|\le|\tilde{u}(t,h(t))-\tilde{u}(t,x_2)|\le M|h(t)-x_2|
 \le 2M|x_1-x_2|.\]

In conclusion, the estimate \eqref{2.8} always holds provided that $0<T\le \min\{1,\frac{1}{2L_1}\}$. Define
 \bess
\mathbb{Y}_{u_0}^T=\{\phi\in C(\Pi_T):\phi(0,x)=u_0(x),\ 0\le\phi\le A,\ |\phi(t,x)-\phi(t,y)|\le 2M|x-y|\}.
 \eess
Obviously, $\mathbb Y_{u_0}^T$ is complete with the metric $d(\phi_1,\phi_2)=\sup_{\Pi_T}|\phi_1-\phi_2|$. For any given $u\in \mathbb Y_{u_0}^T$, we extend $u$ to  $[T,1]\times[-2h_0,2h_0]$ by setting $u(t,x)=u(T,x)$. Then $u\in\mathbb X_{u_0}^1\cap C^{1-}_{x}(\Pi_1)$. Define a mapping $\Gamma$ by
 \[\Gamma(u)=\tilde{u}. \]
The above discussions show that $\Gamma$ maps $\mathbb Y_{u_0}^T$ into itself.

{\color{blue}
{\it Step 3:} We shall show that $\Gamma$ is a contraction mapping in $\mathbb Y_{u_0}^T$ for $0<T\ll1$. Let $(v_i,g_i,h_i)$ be the unique solution of \eqref{2.1} with $u=u_i$, $i=1,2$, and define $t^i_x$ by the manner \eqref{2.3} with $(g,h)=(g_i,h_i)$. Let $\tilde u_i$ be the unique solution of \eqref{2.4} with $t_x=t^i_x$, $v=v_i$ and $T_0=T$. Then
 $$\tilde{u}_i(t,x)=\tilde{u}_i(t^i_x,x)+\int_{t^i_x}^t f_1(s,x,\tilde{u}_i,v_i){\rm d}s
 \ \ \ \mbox{for}\ \ x\in[g_i(T), h_i(T)].$$
Set
 \[U=u_1-u_2, \ \ \widetilde U=\tilde{u}_1-\tilde{u}_2, \ \ h=h_1-h_2, \ \ g=g_1-g_2,\ \  \Omega_T=D^T_{g_1,h_1}\cup D^T_{g_2,h_2}.\]

The following arguments are inspired by those of \cite{DWZ19, LWjmaa}. Make the zero extensions of $\tilde{u}_i$ and $v_i$ in $\big([0,T]\times\mathbb R\big)\setminus D^T_{g_i, h_i}$. Fix $(t,x)\in\Omega_T$, we now estimate $|\widetilde U(t,x)|$ in all the possible cases.

\sk\underline{Case 1}: $x\in(g_1(t),h_1(t))\setminus(g_2(t),h_2(t))$. In such case, either $g_1(t)<x\le g_2(t)$ or $h_2(t)\le x<h_1(t)$, and $\tilde{u}_1(t^1_x,x)=0$, $\tilde{u}_2(t,x)=0$. Thus we have
 \bess
 |\widetilde U(t,x)|=|\tilde{u}_1(t,x)|=\left|\int_{t_x^1}^tf_1(s,x,\tilde{u}_1,v_1){\rm d}s\right|\le C_1|t-t^1_x|,\eess
where $C_1=\|f_1\|_{L^\yy(\Pi_1\times(0,A)\times(0,B))}$.

When $h_2(t)\le x<h_1(t)$, then $0<t^1_x<t$ and $h_1(t)>h_1(t^1_x)=x\ge h_2(t)$. Therefore,
 \bess
 |\widetilde U(t,x)|&\le& C_1|t-t^1_x|=C_1|h_1^{-1}(h_1(t))-h_1^{-1}(h_1(t^1_x))|\\
 &\le& C_1\|(h^{-1}_1)'\|_{L^{\yy}([0,T])}|h_1(t)-h_1(t_x^1)|\\
 &\le& C_1\sigma^{-1}|h_1(t)-h_1(t_x^1)|\\
 &\le&C_1\sigma^{-1}|h_1(t)-h_2(t)|\\
 &\le& C_1\sigma^{-1}\|h\|_{C([0,T])},\eess
where $\sigma>0$ is determined by \eqref{2.5}. When $g_1(t)<x\le g_2(t)$, we can analogously obtain
 \bess
 |\widetilde U(t,x)|=|\tilde{u}_1(t,x)|\le& C_1\sigma^{-1}\|g\|_{C([0,T])}.\eess

\underline{Case 2}: $x\in(g_2(t),h_2(t))\setminus(g_1(t),h_1(t))$. Similar to Case 1 we have
 \bess
 |\widetilde U(t,x)|=|\tilde{u}_2(t,x)|\le C_1\sigma^{-1}\|g,\,h\|_{C([0,T])}.\eess

\sk\underline{Case 3}: $x\in(g_1(t), h_1(t))\cap(g_2(t),h_2(t))$. If $x\in [-h_0, h_0]$, then $t_x^1=t_x^2=0$ and $\tilde u_1(t^1_x,x)=\tilde u_2(t^2_x,x)=\tilde u_0(x)$. Hence
 \begin{align*}
|\widetilde U(t,x)|\le&\int_0^t|f_1(s,x,\tilde{u}_1,v_1)-f_1(s,x,\tilde{u}_2,v_2)|{\rm d}s
\le TL_1\left(\|\widetilde U\|_{C(\up{\Omega}_T)}+\|v_1-v_2\|_{C(\up{\Omega}_T)}\right).
 \end{align*}

If $x\in(g_1(t), h_1(t))\cap(g_2(t),h_2(t))\setminus[-h_0, h_0]$, we have $t_x^1>0, t_x^2>0$, $\tilde u_1(t^1_x,x)=\tilde u_2(t^2_x,x)=0$. Without loss of generality we assume $x>h_0$ and $t_x^2>t_x^1>0$. Then $h_1(t_x^2)>h_1(t_x^1)=x=h_2(t_x^2)$, $x\in(g_1(s), h_1(s))\cap(g_2(s),h_2(s))$ for all $t_x^2<s\le t$ and $x\in(g_1(t_x^2), h_1(t_x^2))\setminus(g_2(t_x^2), h_2(t_x^2))$. Hence,
 \[|\widetilde U(t_x^2,x)|=|\tilde{u}_1(t_x^2,x)|\le C_1\sigma^{-1}
 \|g,\,h\|_{C([0,T])}\]
by the conclusion of Case 1. Integrating the differential equation of $\tilde u_i$ from $t_x^2$ to $t$ we obtain
 \bess
 \tilde u_1(t,x)&=&\tilde{u}_1(t_x^2,x)+\int_{t_x^2}^t f_1(s,x,\tilde{u}_1,v_1){\rm d}s,\\
 \tilde u_2(t,x)&=&\int_{t_x^2}^t f_1(s,x,\tilde{u}_2,v_2){\rm d}s.
 \eess
It follows that
 \bess
 |\widetilde U(t,x)|&\le&\tilde{u}_1(t_x^2,x)+\int_{t_x^2}^t|f_1(s,x,\tilde{u}_1,v_1)-f_1(s,x,\tilde{u}_2,v_2)|{\rm d}s\\
 &\le&|\widetilde U(t_x^2,x)|+\int_0^t|f_1(s,x,\tilde{u}_1,v_1)-f_1(s,x,\tilde{u}_2,v_2)|{\rm d}s\\
&\le&C_1\sigma^{-1}\|g,\,h\|_{C([0,T])}+TL_1\left(\|\widetilde U\|_{C(\up{\Omega}_T)}+\|v_1-v_2\|_{C(\up{\Omega}_T)}\right).
\eess

In conclusion,
 \begin{align}\label{2.9}
 |\widetilde U(t,x)|\le C_1\sigma^{-1}\|g,\,h\|_{C([0,T])}+ TL_1\left(\|\widetilde U\|_{C(\up{\Omega}_T)}+\|v_1-v_2\|_{C(\up{\Omega}_T)}\right).
 \end{align}

We will show in the following Step 4 that if $0<T\ll 1$ then there exists positive constant $C$ such that
\be\label{2.10}
\|g,\,h\|_{C^1([0,T])}\le C\|U\|_{C(\Pi_T)},\ ~\ ~
\|v_1-v_2\|_{C(\bar \Omega_T)}\le C\|U\|_{C(\Pi_T)}.
 \ee
Once this is done, notice that $g(0)=h(0)=0$, the first inequality of \eqref{2.10} implies
 \[\|g,\,h\|_{C([0,T])}\le T\|g,\,h\|_{C^1([0,T])}\le TC\|U\|_{C(\Pi_T)}.\]
Then combing with \eqref{2.9}, we have
 \bess
 \|\widetilde U\|_{C(\Pi_T)}\le \frac 13\|U\|_{C(\Pi_T)} \ \ \mbox{if}\ \ 0<T\ll 1.
\eess}
This demonstrate that $\Gamma$ is a contraction mapping in $\mathbb Y_{u_0}^T$. Thus, $\Gamma$
has a unique fixed point $u$ in $\mathbb Y_{u_0}^T$. Let $(v, g, h)$ be the unique solution of \eqref{2.1} with such $u$. Then $(u, v, g, h)$
is a solution of \eqref{1.1} and it is the unique one provided $u\in\mathbb Y_{u_0}^T$.
Moreover, we can see that $u\in C^{1,1-}(\overline D^T_{g,h})$ and $v\in W_p^{1,2}(D^T_{g,h})$. Thus $v\in C^{1+\alpha,\frac{1+\alpha}{2}}(D^T_{g,h})$ by the embedding theorem as $p>3$.

{\it Step 4:} Proof of \eqref{2.10}.
Its proof is similar to that of \cite[Theorem 2.1: Step 4]{LHWdcdsb} and \cite[Theorem 2.1: Step 3]{WZjde18}.
Before our statement, some preparations are needed. Let
 \bess
 &x_i(t,y)=\dd\frac{1}{2}[(h_i(t)-g_i(t))y+h_i(t)+g_i(t)],&\\[2mm] &\xi_i(t)=\dd\frac{2}{h_i(t)-g_i(t)},\ \ \zeta_i(t,y)=\frac{h'_i(t)+g'_i(t)}{h_i(t)-g_i(t)}
 +\frac{h_i'(t)-g_i'(t)}{h_i(t)-g_i(t)}y,&\\[2mm]
 &w_i(t,y)=u_i(t,x_i(t,y)), \ \ z_i(t,y)=v_i(t,x_i(t,y)), &
 \eess
and
 \[f^i_2(t,y)=f_2(t,x_i(t,y),u(t,x_i(t,y)),v(t,x_i(t,y)))\]
for $i=1,2$. Then,
 \bess
\begin{cases}
z_{i,t}=d\xi^2_iz_{i,yy}+\zeta_iz_{i,y}+f^i_2(t,y),&0<t\le T,\ |y|<1,\\[1mm]
z_i(t,\pm1)=0,&0\le t\le T,\\[1mm]
z_i(0,y)=v_0(h_0y)=:z_{i,0}(y),&|y|\le1.
\end{cases}
\eess
Recall \eqref{2.2}, it follows that
\bes\label{2.11}
\|\xi_i\|_{L^{\yy}([0,T])}\le\frac{1}{h_0},\ \ \|\zeta_i\|_{L^{\yy}(\Delta_T)}\le\frac{2K}{h_0},\ \
\|f_2^i\|_{L^{\yy}(\Delta_T)}\le C_0.
\ees
And
by the $L^p$ theory $\|z_i\|_{W_p^{1,2}(\Delta_T)}\le C_1'$.
Using the arguments in the proof of \cite[Theorem 1.1]{Wdcds19} we can obtain
\bes[z_i,\,z_{i,y}]_{C^{\frac\alpha 2,\alpha}(\Delta_T)}\le \up{C}, \label{2.12}\ees
where $C_1$ is independent of $T^{-1}$.
 This implies
\bes\label{2.13}
\|z_{i,y}\|_{C(\Delta_T)}\le \|z_{i,0}'(y)\|_{C([-1,1])}+\up{C}T^{\frac\alpha 2}
\le \|z_{i,0}'(y)\|_{C([-1,1])}+\up{C}:=C_2'.
\ees
Thanks to \eqref{2.2} and $v_{i,x}(t,x)=z_{i,y}(t,y)\frac{2}{h(t)-g(t)}$, it yields
\bes
\|v_{i,x}\|_{C(\overline D^T_{g,h})}\le {C'_2}/{h_0}.
\label{2.14}\ees

On the other hand, $z=z_1-z_2$ satisfy
\bes\label{2.15}
\left\{\begin{aligned}
&z_t-d\xi_1^2z_{yy}-\zeta_1z_{y}-a(t,y)z
 =d(\xi_1-\xi_2) z_{2,yy}+(\zeta_1^2-\zeta_2^2)z_{2,y}\\
 &\hspace{43mm}+b(t,y)(w_1-w_2)+c(t,y), &&0<t\le T,~|y|<1,\\
&z(t,\pm 1)=0, &&0\le t\le T,\\
&z(0,y)=0, &&|y|\le 1,
\end{aligned}\right.
 \ees
and $g(t)=g_1(t)-g_2(t)$, $h(t)=h_1(t)-h_2(t)$ satisfy
 \bess
\begin{cases}
g'(t)=-\mu\xi_1(t)z_y(t,-1)-\mu(\xi_1(t)-\xi_2(t))z_{2,y}(t,-1), \ \ 0<t\le T,\\[1mm]
h'(t)=-\beta\xi_1(t)z_y(t,1)-\beta(\xi_1(t)-\xi_2(t))z_{2,y}(t,1),\qquad 0<t\le T,\\[1mm]
g(0)=h(0)=0,
\end{cases}
\eess
where
 \bess
 a(t,y)&=&\int_0^1f^1_{2,v}(t,y,w_1,z_2+(z_1-z_2)\tau){\rm d}\tau,\\[2mm]
 b(t,y)&=&\int_0^1f^2_{2,u}(t,y,w_2+(w_1-w_2)\tau,z_2){\rm d}\tau,\\[2mm]
 c(t,y)&=&f^1_2(t,y,w_1,z_2)-f^2_2(t,y,w_1,z_2).
 \eess
Clearly, $\|a,\,b\|_{L^\infty(\Delta_T)}\le L_1$, $\|c\|_{L^\infty(\Delta_T)}\le L^*_1$.
Due to \eqref{2.11}, \eqref{2.13}, applying the parabolic $L^p$ theory to \eqref{2.15} we can obtain
\bess
\|z\|_{W^{1,2}_p(\Delta_T)}\le C_3\big(\|g,\,h\|_{C^1([0,T])}+\|w_1-w_2\|_{C(\Delta_T)}\big),
\eess
where $C_3$ depends on $h_0$, $\mu$, $\beta$, $A$, $B$ and $K$. The same as \eqref{2.12}, we have
\bes
[z]_{C^{\frac\alpha 2,\alpha}(\Delta_T)}+[z_y]_{C^{\frac\alpha 2,\alpha}(\Delta_T)}
&\le& C_4\big(\|g,\,h\|_{C^1([0,T])}+\|w_1-w_2\|_{C(\Delta_T)}\big),
\lbl{2.16}\ees
where $C_4>0$ is independent of $T^{-1}$. When $(t,y)\in \Delta_T$, we have
\bess |w_1(t,y)-w_2(t,y)|&=&|u_1(t,x_1(t,y))-u_2(t,x_2(t,y))|\\[1mm]
 &\le& |u_1(t,x_1(t,y))-u_2(t,x_1(t,y))|+|u_2(t,x_1(t,y))-u_2(t,x_2(t,y))|\\[1mm]
 &\le&\|U\|_{C(\Pi_T)}+L(u,1)|x_1(t,y)-x_2(t,y)|\\[1mm]
 &\le& C_5\kk(\|U\|_{C(\Pi_T)}+\|g,h\|_{C([0,T])}\rr),
\eess
where $C_5$ depends only on $h_0$ and the Lipschitz
constant $L(u,1)$ of $u$. Therefore,
 \[\|w_1-w_2\|_{C(\Delta_T)}\le C_5\big(\|U\|_{C(\Pi_T)}+\|g,h\|_{C([0,T])}\big).\]
This combined with \eqref {2.16} asserts
\bes
[z]_{C^{\frac\alpha 2,\alpha}(\Delta_T)}+[z_y]_{C^{\frac\alpha 2,\alpha}(\Delta_T)}
\le C_6\kk(\|g,\,h\|_{C^1([0,T])}+\|U\|_{C(\Pi_T)}\rr).
\lbl{2.17}\ees
Notice $z_{y}(0,1)=0$. The above estimate implies
\bes
|z_y(t,1)|_{C([0,T])}\le C_6T^{\frac\alpha 2}\kk(\|g,\,h\|_{C^1([0,T])}
+\|U\|_{C(\Pi_T)}\rr).
\lbl{2.18}\ees
As $h(0)=g(0)=0$, it is easy to see that
\bes
|h(t)|\le T\|h'\|_{C([0,T])},\ \ \ |g(t)|\le T\|g'\|_{C([0,T])}.
\label{2.19} \ees
Making use of \eqref{2.13} and \eqref{2.18} we have
\bess
|h'_1(t)-h'_2(t)|&=&\beta|v_{1,x}(t,h_1(t))-v_{2,x}(t,h_2(t))|\nonumber\\[1.5mm]
&=&\beta\left|\df{2[z_{1,y}(t,1)-z_{2,y}(t,1)]}{h_1(t)-g_1(t)}
+2z_{2,y}(t,1)\df{g(t)-h(t)}{[h_1(t)-g_1(t)][h_2(t)-g_2(t)]}\right|\nonumber\\[1.5mm]
&\le&\beta\df{1}{h_0}|z_y(t,1)|+2\beta|z_{2,y}(t,1)|\df{|h(t)|+|g(t)|}{4h_0^2}\nonumber\\[1.5mm]
&\le& C_7T^{\frac\alpha 2}\kk(\|g,\,h\|_{C^1([0,T])}+\|U\|_{C(\Pi_T)}\rr).
\eess
Therefore, by use of \eqref{2.19},
\bess\|h'\|_{C([0,T])}\le C_8T^{\frac\alpha 2}\kk(\|g',\,h'\|_{C([0,T])}+\|U\|_{C(\Pi_T)}\rr).\eess
Similarly, we have
\[\|g'\|_{C([0,T])}\le C'_8T^{\frac\alpha 2}\kk(\|g',\,h'\|_{C([0,T])}+\|U\|_{C(\Pi_T)}\rr).\]
Consequently, $\|g',\,h'\|_{C([0,T])}\le C_9\|U\|_{C(\Pi_T)}$ provided $T$ small enough. Recalling \eqref{2.19} we get {\color{blue}the first inequality of \eqref{2.10}}:
\bes
\|g,\,h\|_{C^1([0,T])}\le C'_9\|U\|_{C(\Pi_T)}.
\label{2.20} \ees
Moreover, as $z(0,y)=0$, we have
 \[|z(t,y)|=|z(t,y)-z(0,y)|\le t^{\frac\alpha 2}[z]_{C^{\frac\alpha 2, \alpha}(\Delta_T)}, \ \ \forall \ (t,y)\in\Delta_T.\]
This combined with \eqref{2.17} allows us to derive
\bes
\|z\|_{C(\Delta_T)}\le T^{\frac\alpha 2}[z]_{C^{\frac\alpha 2, \alpha}(\Delta_T)}\le C_6T^{\frac\alpha 2}\kk(\|g,\,h\|_{C^1([0,T])}+\|U\|_{C(\Pi_T)}\rr).
\lbl{2.21}\ees

Now we estimate $\|v_1-v_2\|_{C(\bar \Omega_T)}$. Fix $(t,x)\in\overline\Omega_T$ let
\[y_i(t,x)=\frac{2x-g_i(t)-h_i(t)}{h_i(t)-g_i(t)}.\]

{\bf Case 1:} $x\in[g_1(t),h_1(t)]\cap[g_2(t),h_2(t)]$. Using \eqref{2.13}, \eqref{2.20} and \eqref{2.21} respectively, we have
\bes|v_1(t,x)-v_2(t,x)|&=&|z_1(t,y_1)-z_2(t,y_2)|\nonumber\\[1.5mm]
&\le&|z_1(t,y_1)-z_2(t,y_1)|+|z_2(t,y_1)-z_2(t,y_2)|\nonumber\\[1.5mm]
&\le&\|z\|_{C(\Delta_T)}+\|z_{2,y}\|_{C(\Delta_T)}|y_1-y_2|\nonumber\\[1.5mm]
&\le&\|z\|_{C(\Delta_T)}+\frac{2}{h_0}\|z_{2,y}\|_{C(\Delta_T)}\|g,h\|_{C([0,T])}\nonumber\\[1.5mm]
&\le& C_6T^{\frac\alpha 2}\big(\|g,\,h\|_{C^1([0,T])}+\|U\|_{C(\Pi_T)}\big)+\frac{2C'_2}{h_0}\|g,h\|_{C([0,T])}\nonumber\\[.5mm]
&\le& C_{10}\|U\|_{C(\Pi_T)}.\label{2.22}
\ees

{\bf Case 2:} $x\in[g_1(t),h_1(t)]\setminus[g_2(t),h_2(t)]$. In this case $v_2(t,x)=0$. Without loss of generality, we may think of $x\in[g_1(t),g_2(t))$ and $g_2(t)\le h_1(t)$. Take advantage of \eqref{2.14} and \eqref{2.22}, it yields
\bess
|v_1(t,x)-v_2(t,x)|&=&|v_1(t,x)-v_2(t,g_2(t))|\\[1mm]
&\le&|v_1(t,x)-v_1(t,g_2(t))|+|v_1(t,g_2(t))-v_2(t,g_2(t))|\\[1mm]
&\le&\|v_{1,x}\|_{C(\overline D^T_{g_1,h_1})}|g_1(t)-g_2(t)|+C_{10}\|U\|_{C(\Pi_T)}\\[1mm]
&\le& C_{11}\|U\|_{C(\Pi_T)}.
\eess

{\bf Case 3:}  $x\in[g_2(t),h_2(t)]\setminus[g_1(t),h_1(t)]$. Similar to Case 2, we still have \[|v_1(t,x)-v_2(t,x)|\le C_{12}\|U\|_{C(\Pi_T)}.\]

In conclusion,
 \[\|v_1-v_2\|_{C(\bar{ \Omega}_T)}\le C\|U\|_{C(\Pi_T)}\]
if $0<T\ll 1$. The estimate \eqref{2.10} is proved.

{\it Step 5:} The uniqueness. Let $(\tilde u, \tilde v, \tilde g, \tilde h)$ be any solution of \eqref{1.1}. It is easy to see from Step 2 that $\tilde u\in\mathbb Y_{u_0}^T$ if $0<T\ll 1$. Thus
$(\tilde u, \tilde v, \tilde g, \tilde h)=(u, v, g, h)$ and the proof is complete.

\vskip 4pt
\begin{proof}[\bf Proof of Theorem \ref{th2.1}]  Clearly, $\phi(t)>0$, $\psi(t)>0$ for all $t\ge 0$. Let $T_*$ be the maximal existence time of $(u,v,g,h)$. For any fixed $0<T<T_*$,  applying the comparison principle in the region $D^T_{g,h}$ we have
 \bes\left\{\begin{array}{lll}
 u(t,x)\le\phi(t)\le\dd\max_{[0,T+1]}\phi(t):=M(T) \ \ \ \mbox{on}\ \ \ol D^T_{g,h}, \\[2mm] v(t,x)\le\psi(t)\le\dd\max_{[0,T+1]}\psi(t):=N(T)
  \ \ \ \mbox{on}\ \ \ol D^T_{g,h}.
  \end{array}\right.\label{2.23}\ees

It is not hard to see that $v_x(t,h(t))<0$ by the Hopf boundary lemma for $0<t<T$, which yields  $h'(t)>0$.
 Set
 $$A=\dd\sup_{[0,M(T)]\times[0,N(T)]}f_2(u,v).$$
Define a comparison function by
 $$w(t,x)=N(T)\left[2K(h(t)-x)-K^2(h(t)-x)^2\right]$$
for some appropriate positive constant $K>1/h_0$ over the region
 $$\Omega_T=\{(t,x):\ 0<t<T,\ h(t)-1/K<x<h(t)\}.$$

First of all, one can easily compute that, for any $(t,x)\in\Omega_T$,
  \[w_t=2N(T)K[1-K(h(t)-x)]h'(t)\geq0,\ \
  -w_{xx}=2N(T)K^2.\]
It follows that, if $K^2\geq \frac{A}{2 dN(T)}$, then
 \[w_t-d w_{xx}\geq 2dN(T)K^2\ge A\geq f_2(t,x,u)=v_t-dv_{xx}\ \ \ {\rm in}\ \ \Omega_T.\]
It is clear that
  $$w(t,h(t)-{K}^{-1})=N(T)\geq v(t,h(t)-{K}^{-1}),\ \ w(t,h(t))=0=v(t,h(t))$$
for all $0<t<T$. Taking advantage of
  \bess
  v_0(x)&=&-\dd\int^{h_0}_xv_0'(y)dy\leq -\min_{[0,h_0]}v_0'(x)(h_0-x),\ x\in[h_0-K^{-1},h_0],\\[1.2mm]
 w(0,x)&\geq& N(T)K(h_0-x),\ \
  x\in[h_0-K^{-1},h_0],\eess
we have that if
 \[N(T)K\geq-\min_{[0,h_0]}v_0'(x)\]
then
  $$v_0(x)\leq w(0,x) \ \ \mbox{in} \ \, [h_0-K^{-1},h_0].$$
Applying the maximum principle to $w-v$ over $\Omega_T$ we deduce $w\ge v$ in $\Omega_T$. It then leads to $v_x(t,h(t))\geq w_x(t,h(t))=-2N(T)K$. Thus we have
 \bes
 h'(t)=-\beta v_x(t,h(t))\le 2\beta N(T)K=2\beta\max\left\{\frac 2{h_0}, \ \sqrt{\frac{AN(T)}{2d}}, \
 -\min_{[0,h_0]}v_0'(x)\right\}\label{2.24}\ees
for all $0\le t\le T$. Similarly,
  \bes
  g'(t)\ge -2\mu\max\left\{\frac 2{h_0}, \ \sqrt{\frac{AN(T)}{2d}}, \
 \max_{[-h_0, 0]}v_0'(x)\right\}, \ \ \forall \ 0\le t\le T.\label{2.25}\ees

Recalling the estimates \eqref{2.23}-\eqref{2.25} and using a similar method to the proof of \cite[Theorem 1.2]{Wdcds19} we have  $T_*=\yy$. \end{proof}

\vskip 4pt
\begin{proof}[\bf Proof of Theorem \ref{th2.2}] It is easy to see that
 \[u(t,x)\le\dd\max_{[-h_0,h_0]} u_0+k_0:=\eta, \ \ v(t,x)\dd\le \max_{[-h_0,h_0]} v_0+\Theta(\eta).\]
The remaining proof is the same as that of Theorem \ref{th2.1}. \end{proof}

\end{document}